\documentclass[12pt,thmsa]{amsart}

\def\bbr{{\Bbb R}}

\def\vol{{\hbox{\rm vol}}}

\def\part{\partial}
\def\intl{\int\limits}
\def\b{\beta}

\def\Gam{\Gamma}

\def\a{\alpha}

\def\del{\delta}
\def\vp{\varphi}

\def\gam{\gamma}

\def\sig{\sigma}

\def\e{\varepsilon}

\newtheorem{theorem}{Theorem}[section]
\newtheorem{lemma}[theorem]{Lemma}

\newtheorem{corollary}[theorem]{Corollary}

\theoremstyle{remark}
\newtheorem{remark}[theorem]{Remark}

\newtheorem{example}[theorem]{Example}

\numberwithin{equation}{section}

\newcommand{\be}{\begin{equation}}
\newcommand{\ee}{\end{equation}}

\newcommand{\bea}{\begin{eqnarray}}
\newcommand{\eea}{\end{eqnarray}}
\newcommand{\Bea}{\begin{eqnarray*}}
\newcommand{\Eea}{\end{eqnarray*}}

\begin{document}

\title[The generalized Busemann-Petty problem]
{The generalized Busemann-Petty problem with weights
 }

\author{Boris Rubin}
\address{
Department of Mathematics, Louisiana State University, Baton Rouge,
LA, 70803 USA}

\email{borisr@math.lsu.edu}


\subjclass[2000]{Primary 52A38; Secondary 44A12}



\keywords{The generalized Busemann-Petty problem with weights, Radon
transforms, star bodies}

\begin{abstract}
The generalized Busemann-Petty  problem  asks whether
origin-symmetric convex bodies  in $\bbr^n$ with smaller
$i$-dimensional sections necessarily have smaller volume. We study
 the weighted version of this problem corresponding to the physical
 situation when bodies
 are endowed with  mass distribution and the relevant
 sections are measured with attenuation.
\end{abstract}

\maketitle

\section{Introduction}

\setcounter{equation}{0}

Let $G_{n,i}$  be the Grassmann manifold of $i$-dimensional linear
subspaces  of $\bbr^n$, and let $\vol_i ( \cdot)$ denote the
$i$-dimensional volume function, $1 \le i \le n$. Is it true that
for origin-symmetric convex bodies $K$ and $L$ in $\bbr^n$, the
inequality \be \vol_i(K \cap \xi) \le \vol_i(L \cap \xi) \quad
\forall \xi \in G_{n,i} \ee implies  \be \vol_n(K)  \le \vol_n(L)
\quad \text{\rm ?} \ee This  question is known as {\it the
generalized Busemann-Petty problem}. For $i=n-1$,  the problem was
posed by Busemann and Petty \cite{BP} in 1956. It has a long
history, and the answer is affirmative if and only if $n \le 4$; see
\cite{G}, \cite{K3}, \cite{R2}. For the generalized Busemann-Petty
problem the following statements are known. If $i=2, n= 4$, an
affirmative answer follows from that in the case $i=n-1$. If $3<i
\le n-1$, the
 negative answer was given by Bourgain and  Zhang \cite{BZ}; see also
 \cite{K3}, \cite{RZ}.
  For  the special case, when
  $K$ is a body of revolution, the answer for $i=2$ and $3$
  is affirmative \cite{GZ}, \cite{Z2}, \cite{RZ}.
  The  case, when $K$ is an arbitrary origin-symmetric convex
body and $i=2$ and $3$, is still open.

In a recent paper \cite{Zv}, Zvavitch  considered the Busemann-Petty
problem ($i=n-1$) in a more general setting, when volumes under
consideration are evaluated with respect to general measures
satisfying certain conditions; see also \cite{Zv1} where the case of
the Gaussian measure was considered. Motivated by these papers, we
extend the results from \cite{Zv} to sections of arbitrary dimension
 $1\le i\le n-1$ and study a weighted version of the generalized
 Busemann-Petty problem. Our approach is new in the sense that it
 relies on elementary
properties of Radon transforms on the sphere and does not invoke the
Fourier transform techniques as in \cite{Zv}.  Main results are
presented by Theorems \ref{th1} and \ref{th3}. Diverse geometric
inequalities that follow from those theorems are exhibited in
Section 4.

The generalized Busemann-Petty problem with weights can be given a
physical meaning, when bodies under consideration
 are endowed with  mass distribution and the relevant
 sections are measured with  inevitable attenuation.

I would like to thank Prof. Alex Koldobsky for useful discussions.

\section{Preliminaries}

\setcounter{equation}{0}
 We use the following notation:
 $S^{n-1}$ is the unit sphere in $\bbr^n$; $\, \sig_{n-1}=
2\pi^{n/2}/\Gam (n/2) $ is the area of $S^{n-1}$;  $ \, e_1, e_2,
\ldots , e_n$ denote the coordinate unit vectors.  In the following
$SO(n)$ is the special orthogonal group of $\bbr^n$; $ \, SO(n-1)$
stands for
 the subgroup of $SO(n)$ preserving $e_n$. If $i$ is an integer, $1\le i\le
n-1$, then  $G_{n,i}$ denotes the Grassmann manifold of
$i$-dimensional linear subspaces
 of $\Bbb R^n$. For $\gam \in SO(n)$, and
$\xi \in G_{n,i}$, we denote by  $d\gam $ and $d\xi$ the
corresponding $SO(n)$-invariant measures with total mass $1$.

For continuous functions $f(\theta)$ on $S^{n-1}$ and $\varphi
(\xi)$ on $G_{n,i}$, the totally geodesic Radon transform $R_i f$
 and its
dual  $R_i^*\varphi$ are defined by \be\label{rts}
 (R_i f)(\xi) = \intl_{S^{n-1}\cap\xi} f(\theta) \, d_\xi \theta, \qquad
  (R_i^* \varphi)(\theta) = \intl_{\xi \ni \theta}  \varphi (\xi)  \, d_\theta \xi,
\ee where $d_\xi \theta$ and $ d_\theta \xi$ denote the induced
 measures on the corresponding manifolds $S^{n-1}\cap\xi$
and $\{\xi\in G_{n,i}: \xi \ni \theta \}$; see \cite{He}, \cite{R1}.
The precise meaning of the second integral is
 \be\label{drt}
 (R_i^*
 \vp)(\theta)=\int_{SO(n-1)}\vp (r_\theta \gam p_0) \,  d\gam, \qquad \theta \in
 S^{n-1},
 \ee
 where $p_0 =\bbr
e_{n-i+1} + \ldots + \bbr e_{n}$ is the coordinate $i$-dimensional
plane and $ r_\theta \in SO(n)$ is a rotation satisfying $r_\theta
e_n =\theta$. The corresponding duality relation
 reads
\be\label{dual} \frac{1}{\sig_{i-1}} \intl_{G_{n,i}} (R_if)(\xi) \vp
(\xi) d \xi = \frac{1}{\sig_{n-1}} \intl_{S^{n-1}} f(\theta) (R_i^*
\vp) (\theta) d\theta \ee and is applicable provided the integral in
either side is finite for $f$ and $\vp$ replaced by $|f|$ and
$|\vp|$, respectively.

The Radon transform $R_i$ and its dual extend as linear bounded
operators from $L^1(S^{n-1})$ to $L^1(G_{n,i})$ and from
$L^1(G_{n,i})$ to $L^1(S^{n-1})$, respectively. Moreover, they can
be defined for finite Borel measures. Specifically, if $\mu$ is such
a measure on $G_{n,i}$, then, according to  (\ref{dual}), $R_i^*
\mu$ is a finite Borel measure on $S^{n-1}$ (i.e., a linear
continuous functional on $C(S^{n-1})$) defined by
 \be\label{mer}(R_i^* \mu, f) =
 \frac{\sig_{n-1}}{\sig_{i-1}}\intl_{G_{n,i}}
(R_i f)(\xi)  \, d\mu (\xi), \quad f \in C(S^{n-1}).\ee For
instance, if $\mu$ is the unit mass on the circle $S^{n-1}\cap
\bbr^i$, then $R_i^* \mu$ assigns to $f$ the integral of $f$ over
this circle multiplied by $\sig_{n-1}/\sig_{i-1}$.

Let $K$ be an origin-symmetric star body  in $\bbr^n$.  The {\it
radial function} of $K$ is defined by
\[
\rho_K (\theta) = \sup \{ \lambda \ge 0: \, \lambda \theta \in K \},
\qquad \theta \in S^{n-1},\] and represents the Euclidean distance
from the origin to the boundary of $K$ in the direction of $\theta$.
If $\xi$ is an $i$-dimensional subspace of $\bbr^n$, $ 1\le i\le n$,
then \be\label{vol} \vol_i(K\cap \xi) = i^{-1}\intl_{S^{n-1}\cap
\xi} \rho_K^i (\theta) \,d_\xi \theta.\ee If $1\le i< n$ this is
just $i^{-1}(R_i \rho_K^i )(\xi)$. The body $K$ is called infinitely
smooth if $\rho_K(\theta)\in C^\infty_{even}(S^{n-1})$.

We will need the following elementary inequality which is a slight
generalization of Lemma 1 from \cite{Zv}.
\begin{lemma}\label {zl} Let $a, b>0$
and suppose that $\a(r)$ and $\b(r)$ are positive continuous
functions on $(0, \max \{a,b\})$ such that $r^{n-i}\a(r)/\b(r)$ is
nondecreasing on $(0, \max \{a,b\})$. Then \bea &&\intl_0^a
r^{n-1}\a(r)\, dr -a^{n-i}\frac{\a (a)}{\b(a)}\intl_0^a
r^{i-1}\b(r)\, dr \nonumber \\ \label{zvn}&\le& \intl_0^b
r^{n-1}\a(r)\, dr-a^{n-i}\frac{\a (a)}{\b(a)}\intl_0^b
r^{i-1}\b(r)\, dr.\eea \end{lemma}\begin{proof} This inequality is
equivalent to
\[a^{n-i}\frac{\a (a)}{\b(a)}\intl_a^b
r^{i-1}\b(r)\, dr\le \intl_a^b r^{n-1}\a(r)\, dr.\] The latter is
obvious by taking into account that $r^{n-i}\a(r)/\b(r)$ is
nondecreasing, no matter $a<b$ or $a>b$.
\end{proof}

\section{Main theorems}

Let $K$ be an origin-symmetric star body in $\bbr^n$
 with the radial function $\rho_K (\theta)$. Given nonnegative
 measurable functions $u$ and $v$ on $\bbr^n$, we
  denote \[ V_u (K\cap \xi)=\intl_{K\cap \xi}
 u(x)\,dx, \qquad V_v (K)=\intl_{K}
 v(x)\,dx,\] provided these integrals are well defined.
 The functions $u$ and $v$ can be given a physical meaning
 to be the attenuated mass distribution
 and the true mass distribution, respectively.
In polar coordinates we have \[ V_u (K\cap \xi) \!=
\!\intl_{S^{n-1}\cap \xi} \!   \!\!d\theta  \! \!
 \intl_0^{\rho_K (\theta)}  \!\!r^{i-1}u(r \theta)\,dr, \quad V_v (K) \!=
 \!\intl_{S^{n-1}}  \! \!d\theta  \! \!
 \intl_0^{\rho_K (\theta)}  \!\!r^{n-1}v(r \theta)\,dr.\]
The second integral  is finite for any locally integrable function
$v$. The first one is represented as the Radon transform \be
\label{rtr} V_u (K\cap \xi)=(R_i b_K)(\xi), \qquad b_K
(\theta)=\intl_0^{\rho_K (\theta)}r^{i-1}u(r
 \theta)\,dr.\ee
It is finite (at least for almost all $\xi \in G_{n,i}$) if
$|x|^{i-n}u(x)$ is locally integrable.
 This follows from duality (\ref{dual}), according
to which (set $\vp \equiv 1$)\bea \frac{\sig_{n-1}}{\sig_{i-1}}
\intl_{G_{n,i}}V_u (K\cap \xi)\, d\xi
&=&\frac{\sig_{n-1}}{\sig_{i-1}} \intl_{G_{n,i}}(R_i b_K)(\xi)
= \intl_{S^{n-1}}b_K (\theta)\,d\theta\nonumber\\
&=& \intl_{S^{n-1}}d\theta \intl_0^{\rho_K (\theta)}r^{i-1}u(r
 \theta)\,dr= \intl_{K}|x|^{i-n}u(x)\, dx.\nonumber\eea

For technical reasons  we impose some more restrictions on $u$ and
$v$ and consider a class of weights satisfying the
 following conditions:

{\rm (a)} $u(x)$ is an even function which is  positive and
continuous for $x \in \bbr^n \setminus \{0\}$ and such that
 $|x|^{i-n}u(x)$ is locally integrable;

{\rm (b)} $v(x)$ is a nonnegative, even, locally integrable
function, and the function $v_\theta(r)=v(r\theta)$ is continuous in
$r>0$ for almost all $\theta \in S^{n-1}$;

{\rm (c) ({\bf the comparison condition})} The function  $ \;
a_\theta(r)=r^{n-i}\,\frac{v(r \theta)}{u(r \theta)}$ is
nondecreasing  for almost all $\theta \in S^{n-1}$.

The conditions (a)-(c) look pretty sophisticated but they allow us
to consider  weights $v$ which are discontinuous on the unit sphere;
see Example 4.4. The comparison condition (c) restricts our  class
of admissible weights,  and the case when (c) fails remains open.
However, this condition has a certain physical meaning: if
attenuation is too strong,  we cannot retrieve desired information
from measurements.

   Given  a symmetric star body $K$ in $\bbr^n$, we introduce a
   {\it comparison function}
 \be\label{ak} a_K(\theta) \equiv a_\theta(\rho_K (\theta))=
 \rho_K^{n-i} (\theta)
 \,  \frac{v(\rho_K (\theta)\, \theta)}{u(\rho_K (\theta)\, \theta)}. \ee

\begin{theorem}\label{th1} Let $2\le i\le n-1$ and suppose that
$u$ and $v$ satisfy the conditions (a)-(c) above.
 If the comparison function $a_K(\theta)$ is represented by the dual Radon transform of a
 positive measure $\mu$ on
 $G_{n,i}$, i.e., $a_K=R_i^* \mu$, then for any symmetric
 star body  $L$ in $\bbr^n$, satisfying \be\label{per} \intl_{K\cap \xi}
 u(x)\,dx \le \intl_{L\cap \xi}
 u(x)\,dx, \qquad \forall \xi \in G_{n,i},\ee
 we have
 \be\label{vt} \intl_{K}
 v(x)\,dx \le \intl_{L}
 v(x)\,dx.\ee
\end{theorem}
A few words are in order  on how one should  interpret
 the key equality $a_K=R_i^* \mu$. Note that by (a) and (b), the functions $b_K$
and $b_L$ are continuous, and $a_K \in L^1(S^{n-1})$. On the other
hand, $R_i^* \mu$ is a measure; see definition (\ref{mer}). The
equality $a_K=R_i^* \mu$ means that $\int_{S^{n-1}}a_K(\theta)
f(\theta)\, d\theta=(R_i^* \mu, f)$ for any $f \in C(S^{n-1})$ or
$R_i^* \mu$ is an absolutely continuous measure (with respect to the
Lebesgue measure on $S^{n-1}$) with density $a_K$.

{\it Proof of Theorem }\ref{th1}. The result is an immediate
consequence of the following inequalities: \be\label{eeq1}
\intl_{S^{n-1}} a_K(\theta)\,b_K (\theta)\,d\theta
\le\intl_{S^{n-1}} a_K(\theta)\,b_L (\theta)\,d\theta.\ee
\be\label{eeq2}V_v(K) -\intl_{S^{n-1}}
  a_K (\theta) b_K (\theta)\,d\theta\le V_v(L)
  - \intl_{S^{n-1}} a_K (\theta) b_L (\theta)\,d\theta,\ee
in which $a_K(\theta), \; b_K (\theta)$ and $b_L (\theta)$ are
defined by (\ref{ak}) and (\ref{rtr}). The
 inequality (\ref{eeq1}) can be easily obtained if we write
(\ref{per}) as $(R_i b_K)(\xi) \le(R_i b_L)(\xi)$ and make use of
the definition (\ref{mer}):
 \bea &&\intl_{S^{n-1}} a_K(\theta)\,b_K (\theta)\,d\theta
 =(R_i^* \mu, b_K)=\frac{\sig_{n-1}}{\sig_{i-1}}\intl_{G_{n,i}}(R_i b_K)(\xi)\,
d\mu(\xi)\nonumber\\
&&\le \frac{\sig_{n-1}}{\sig_{i-1}}\intl_{G_{n,i}}(R_i b_L)(\xi)\,
d\mu(\xi)=(R_i^* \mu, b_L)=\intl_{S^{n-1}} a_K(\theta)\,b_L
(\theta)\,d\theta.\nonumber\eea The inequality  (\ref{eeq2})  can be
   derived from (\ref{zvn}) if we set $\; a=\rho_K
(\theta)$,  $\;b=\rho_L (\theta)$,  $\;\a(r)=v(r \theta), \;
\b(r)=u(r \theta)$. This gives \bea
 &&\intl_0^{\rho_K (\theta)}r^{n-1}v(r \theta)\,dr-\rho_K^{n-i} (\theta)
 \,  \frac{v(\rho_K (\theta)\, \theta)}{u(\rho_K (\theta)\,
 \theta)}\, \intl_0^{\rho_K (\theta)}r^{i-1}u(r \theta)\,dr
 \nonumber \\&\le& \intl_0^{\rho_L (\theta)}r^{n-1}v(r \theta)\,dr-\rho_K^{n-i} (\theta)
 \,  \frac{v(\rho_K (\theta)\, \theta)}{u(\rho_K (\theta)\,
 \theta)}\, \intl_0^{\rho_L (\theta)}r^{i-1}u(r \theta)\,dr\nonumber
 \eea or
\[\intl_0^{\rho_K (\theta)}r^{n-1}v(r \theta)\,dr-a_K(\theta)\,b_K
(\theta) \le \intl_0^{\rho_L (\theta)}r^{n-1}v(r
\theta)\,dr-a_K(\theta)\,b_L (\theta).\]
  Integrating the latter over $S^{n-1}$, we obtain (\ref{eeq2}).
  $\hfill\hfill\square$

\begin{remark} 1.  We did not include the case $i=1$ in Theorem
\ref{th1} because in this case the implication
(\ref{per})$\Rightarrow $(\ref{vt}) is true for any nonnegative $u$
and $v$ satisfying the condition (a) and (b) above.
\end{remark}

The next theorem shows that the assumption $a_K=R_i^* \mu,\; \mu>0$,
in Theorem \ref{th1} is crucial. Namely, if it fails, then there
exist origin-symmetric convex bodies $K$ and $L$ such that  $V_u
(K\cap \xi)\le  V_u (L\cap \xi)$ for all  $\xi \in G_{n,i}$, but $
V_v (K)> V_v (L)$. More precisely, the following statement holds.

\begin{theorem}\label{th3} Let $u$ and $v$ satisfy the conditions
(a)-(c) above. Suppose also that $v$ is  positive and both functions
are infinitely differentiable away from the origin.
  Given  an infinitely smooth origin-symmetric convex body $L \subset \bbr^n$
  with positive curvature, let
 \be\label{al} a_L(\theta) \equiv a_\theta(\rho_L (\theta))=
 \rho_L^{n-i} (\theta)
 \,  \frac{v(\rho_L (\theta)\, \theta)}{u(\rho_L (\theta)\, \theta)} \ee
 be represented by the dual Radon transform $R_i^* \vp$ of a function
 $\vp \in C^\infty(G_{n,i})$ which
 is negative for some  $\xi \in G_{n,i}$.
 Then there is  a convex
 symmetric  body $K$ in $\bbr^n$
 such that  \be\label{perr} \intl_{K\cap \xi}
 u(x)\,dx \le \intl_{L\cap \xi}
 u(x)\,dx, \qquad \forall \xi \in G_{n,i},\ee
 but
 \be\label{vtt} \intl_{K}
 v(x)\,dx > \intl_{L}
 v(x)\,dx.\ee
\end{theorem}
\begin{proof} We start with some comments that might be useful for
understanding the essence of the matter. Since the mapping $R_i^* :
C^\infty(G_{n,i}) \to C^\infty(S^{n-1})$ is ``onto", the function
$a_L(\theta)$ is represented as $R_i^* \vp$ for some $\vp \in
C^\infty(G_{n,i}) $ automatically. Such a function $\vp$ is not
unique for $1<i<n-1$, because $R_i^*$ is non-injective in this case.
The theorem actually assumes that there as at least one
representative of the class $\{\vp +\ker (R_i^*)\}$ which is
negative somewhere on $G_{n,i}$.

 As in the previous theorem,
the result will follow if define $K$ satisfying the following
inequalities: \be\label{eeeq1} \intl_{S^{n-1}} a_L(\theta)\,b_K
(\theta)\,d\theta
> \intl_{S^{n-1}} a_L(\theta)\,b_L (\theta)\,d\theta,\ee
\be\label{eeeq2}V_v(K) -\intl_{S^{n-1}}
  a_L (\theta) b_K (\theta)\,d\theta \ge V_v(L)
  - \intl_{S^{n-1}} a_L (\theta) b_L (\theta)\,d\theta.\ee
The body $K$ can be defined as follows. Since $\vp$
 is smooth, then there exist  $\del >0$ and $\theta_0 \in S^{n-1}$
  such that $\vp(\xi)$  is negative for all $\xi$ in the
  open domain  $\Omega_\del =\{\xi \in G_{n,i} :
d(S^{n-1}\cap\xi, \theta_0)<\del
  \}$, $d(\cdot, \cdot)$ being the  geodesic distance on $S^{n-1}$. Consider
  the spherical cap $B=\{\theta : d(\theta, \theta_0)<\del
  \}$, and  let $B'$ denote the symmetric cap centered at
  $-\theta_0$. Choose a non-positive function $g \in
C^\infty_{even}(S^{n-1})$, $g \not\equiv 0$,
 supported by $B\cup B'$
  Then $g_1=R_i g$ is a $C^\infty$ negative function supported by
 $\Omega_\del$, and by duality (\ref{dual}) we have
 \be\label{neg} \intl_{S^{n-1}}
  a_L \, g =  \intl_{S^{n-1}} g \, R_i^* \vp=
  \frac{\sig_{n-1}}{\sig_{i-1}} \intl_{G_{n,i}} g_1 \, \vp >0.\ee
Now we define  an  origin-symmetric convex body $K$
 so that \be \label{dk} b_K (\theta)= b_L (\theta)+\e
g(\theta),\ee assuming $\e>0$ sufficiently small (the proof of
validity of this definition almost coincides with that of
Proposition 2 in \cite{Zv}). Multiplying (\ref{dk}) by $a_L$ and
integrating over $S^{n-1}$, we get
\[  \intl_{S^{n-1}}  a_L \, b_K =\intl_{S^{n-1}}  a_L \, b_L +
\e\intl_{S^{n-1}} a_L \, g.\] Owing to (\ref{neg}), this gives
(\ref{eeeq1}). The proof of  (\ref{eeeq2}) is similar to that of
(\ref{eeq2}) in Theorem \ref{th1} and relies on the inequality
(\ref{zvn}) in which one should set $a=\rho_L (\theta), \; b=\rho_K
(\theta), \; \a(r)=v(r \theta), \; \b(r)=u(r \theta)$.
\end{proof}

\section{Corollaries and partial results}

Theorems \ref{th1} and \ref{th3} give rise to a  series of
statements. Some of them are new and others were obtained before in
a more complicated way.  Below we present a few examples.

\subsection{The case of equal weights}

Let $u$  be a positive even functions on $\bbr^n$ which is
continuous away from the origin and  $|x|^{i-n}u(x)$ is locally
integrable.  Suppose the weights in Theorems \ref{th1} and \ref{th3}
are equal, i.e., $v\equiv u$. Then $a_K (\theta)=\rho_K^{n-i}
(\theta)$ and we have the following statement.
\begin{corollary}\label{urv} $\qquad $$\qquad $
$\qquad $$\qquad $$\qquad $$\qquad $$\qquad $

{\rm (i)} If $\rho_K^{n-i}=R_i^* \mu$ where $\mu$ is a positive
measure on the Grassmanninan $G_{n,i}$\footnote{Origin-symmetric
star bodies with this property were called in \cite{Z1}
$i$-intersection bodies. See  \cite{Lu} and \cite{GLW} for
$i=n-1$.}, then for any symmetric star body $L$ in $\bbr^n$,
satisfying \be\label{cr1} V_u (K\cap \xi)\le V_u (L\cap \xi)\quad
\forall \xi \in  G_{n,i} \ee we have $ V_u (K)\le V_u (L)$.

{\rm (ii)} Let $L$ be an infinitely smooth origin-symmetric convex
body in $\bbr^n$ so that $\rho_L^{n-i}=R_i^* \vp$  for some $\vp \in
C^\infty (G_{n,i})$. If $\vp (\xi)<0$ for some $\xi \in  G_{n,i}$ ,
then there is a convex
 symmetric  body $K$ in $\bbr^n$ which obeys (\ref{cr1}) and $ V_u (K)> V_u (L)$.
 \end{corollary}

For $i=n-1$ this statement was proved by A. Zvavitch \cite{Zv} who
used the Fourier transform approach. The key question is what can
one say about  validity of the representation
\be\label{rok}\rho_K^{n-i}=R_i^* \mu,\qquad \mu \ge 0.\ee
 It is known \cite {BZ},
\cite {K2}, \cite {RZ}, that if $i>3$, then  there is an infinitely
smooth origin-symmetric strictly convex body for which (\ref{rok})
fails, and  we are in the situation of the statement (ii) above. In
the special case $i=n-1$ corresponding to the Busemann-Petty problem
with equal weights, this gives a negative answer to this problem for
all $n>4$. If $n=3,4$, the validity of (\ref{rok}) for $i=n-1$ was
proved by different methods in a series of publications; see, e.g.,
\cite{G}, \cite {K3},  \cite {Z2}, \cite {R2}, and references
therein.

The cases $i=2$ and  $i=3$ when $n>4$ are the most difficult. In
these cases the validity of (\ref{rok}) is known only for bodies of
revolution \cite {GZ}, \cite {RZ}. For arbitrary convex bodies the
problem is still open.

\subsection{The case of power weights}

Let $u(x)=|x|^\a, \; v(x)=|x|^\b$. Then the conditions (a)-(c) have
the form \be\label{ab1} 0<\a +i \le \b +n. \ee The function $a_K
(\theta)$ is $\rho_K (\theta)^{\b+n-\a-i}$. Representation of this
function by the dual Radon transform of a positive measure and the
relevant generalization of the Busemann-Petty problem was studied in
\cite{RZ}. By making use of Erdelyi-Kober fractional integrals, it
was proved, that for every $i>3$, there exist an infinitely smooth
origin-symmetric strictly convex body $L$ of revolution for which
the representation $\rho_L (\theta)^{\b+n-\a-i}=R^*_i \mu$ fails to
be true with $\mu>0$. By Theorem \ref{th3}, it follows that if $i>3$
and $0<\a +i \le \b +n$, then there exists a convex symmetric body
$K$ such that \be\label{cr2} \intl_{K\cap \xi}\!
 |x|^\a\,dx \le\! \intl_{L\cap \xi}\!
 |x|^\a\,dx \; \forall \xi \in G_{n,i},\qquad \!\intl_{K}\!
 |x|^\b\,dx > \intl_{L}
 |x|^\b\,dx.\ee
  For $i=2$ and $3$, the representation $\rho_K
(\theta)^{\b+n-\a-i}=R^*_i \mu$,
 $\mu>0$, corresponding  to Theorem \ref{th1}, is known to be true in
 the case $\a +i+1 = \b +n$ \cite{RZ}. We observe  that it is
 also true if  $\a +i =
\b +n$ because in this  case the equality $1=R^*_i \mu$ trivially
holds with $\mu \equiv 1$. More subtle  results in the cases $i=2$
and $3$, covering the whole domain (\ref{ab1}), were obtained for
bodies of revolution; see \cite{RZ} for details. For arbitrary
symmetric convex bodies, the case $\a +i \neq \b +n$ ($i=2,3)$
remains open. The case $\a +i > \b +n$ contradicts (\ref{ab1}) and
is also open because it does not fall into the scope of Theorems
\ref{th1} and \ref{th3} (in this case the condition (c) is not
satisfied).

It is worth exhibiting the particular case $\b=0;\; i=2,3$, when the
implication \be\label{imp} \intl_{K\cap \xi}
 |x|^\a\,dx \le \intl_{L\cap \xi}
 |x|^\a\,dx \quad \forall \xi \in G_{n,i} \;\Longrightarrow \;\vol_n(K)\le \vol_n(L)\ee
holds provided $\a=n-i-1$ and $\a=n-i$. It may fail if $\a<0$ and
the question is open for $0\le \a<n-i \;(\a \neq n-i-1)$ and
$\a>n-i$.

\subsection{More general homogeneous weights}

The case $\a-\b=n-i$ in the previous subsection when $a_K
(\theta)\equiv 1$ deserves special mentioning. In this case, owing
to  Theorem \ref{th1},
 the implication \be\label{per1} \intl_{K\cap \xi}|x|^\a
 \,dx \le \intl_{L\cap \xi}
 |x|^\a\,dx \; \forall \xi \in G_{n,i}\;\Longrightarrow \;
  \intl_{K}|x|^\b\,dx \le \intl_{L}|x|^\b \,dx\ee
is valid {\it for all} symmetric star bodies $K$ and $L$ and all
$0<i<n$. This  observation  can be essentially generalized. One can
ask the following  question: For which more general homogeneous
weights the implication \be\label{gen1} V_u (K\cap \xi)\le V_u
(L\cap \xi) \; \forall \xi \in G_{n,i}\;\Longrightarrow \; V_v
(K)\le V_v (L)\ee is independent of the choice of  symmetric star
bodies $K$ and $L$, i.e., $a_K (\theta)$ is independent of $K$? The
following theorem answers this question.

\begin{theorem}\label{th4}
 Let $u $ and $v$ be homogeneous functions of
degree $\a$ and $\b$, respectively, which satisfy the conditions
(a)-(c) above. Suppose that  $\a-\b=n-i$ and there is a function
$\vp\in L^1(G_{n,i})$ such that \be\label{cnd1}
v(\theta)=u(\theta)\,(R_i^* \vp)(\theta) \ee for almost all $\theta
\in S^{n-1}$. Then the implication (\ref{gen1}) holds for any
symmetric star bodies $K$ and $L$ in $\bbr^n$.
 \end{theorem}
\begin{proof} The  statement is a consequence of Theorem \ref{th1},
because for any symmetric star bodies $K$, \[a_K(\theta)=
 \rho_K^{n-i} (\theta)
 \,  \frac{v(\rho_K (\theta)\, \theta)}{u(\rho_K (\theta)\,
 \theta)}=\rho_K^{n-i} (\theta)
 \,  \frac{\rho_K^{\b}(\theta)\,v(\theta)}{\rho_K^{\a}(\theta)\,u(
 \theta)}= \frac{v(\theta)}{u(
 \theta)}  =(R_i^*\vp)(\theta).\]
\end{proof}

\begin{example}Let us consider the weight functions
$$u(x)=|x|^\a,  \qquad v(x)=|x|^\b\, w_\gam(x),$$ where $$
w_\gam(x)=
(1-x^2_n/|x|^2)^{(\gam+i-n)/2}=(|x'|/|x|)^{\gam+i-n},\quad x' =(x_1,
\ldots, x_{n-1}).$$
 Suppose that
 $\a>-i, \; \a-\b=n-i$, and $\gam >0$. It is known (see Example 2.5 in \cite {R1})
  that $ w_\gam(\theta)=(R^*_i m_\gam)(\theta)$ with $$m_\gam(\xi)=\frac{\sig_{n-2} \,
\Gam((i-1+\gam)/2)}{\pi^{(i-1)/2} \, \sig_{n-i-1} \, \Gam(\gam/2)}\,
\sin^{\gam+i-n} [d(e_{n}, S^{n-1}\cap\xi)],$$
  $d(\cdot,\cdot)$ being
 the geodesic distance on $S^{n-1}$.  By
Theorem \ref{th4},
 for any symmetric star bodies $K$ and $L$ in $\bbr^n$, the
 inequality
\be \intl_{K\cap \xi} |x|^\a
 \,dx \le \intl_{L\cap \xi}
 |x|^\a\,dx \quad \forall \xi \in G_{n,i}\ee
 implies
\be\intl_{K}|x|^\b \,w_\gam(x) \,dx \le \intl_{L}|x|^\b \,w_\gam(x)
\,dx.\ee In particular (set $\a=0$) for any $\gam >0$,
\be\intl_{K}|x'|^{\gam+i-n}|x|^{-\gam}\,dx \le
\intl_{L}|x'|^{\gam+i-n}|x|^{-\gam}\,dx\ee provided $\int_{K\cap
\xi} dx \le \int_{L\cap \xi}dx \quad \forall \xi \in G_{n,i}$.
\end{example}
We conclude this article by laying stress on the question that is of
major importance in  Theorems \ref{th1} and \ref{th3}: Is it
possible to represent the comparison function $a_K (\theta)$ by the
dual Radon transform of a positive measure? This question is
difficult even in the case $i=n-1$ when the corresponding Radon
transform (it is known as  the Minkowski-Funk transform) is actually
self-adjoint and injective. The case $1<i<n-1$ is much more
difficult because the dual Radon transform is non-injective for such
$i$ (it has a nontrivial kernel). These difficulties have been
overcome so far only in some particular cases using the tools
fractional calculus, the Fourier analysis, and known facts from the
theory of Radon transforms.

\end{document}